\documentclass[12pt,reqno]{amsart}
\usepackage{amsmath, amssymb, amsfonts, xcolor}
\usepackage{hyperref}
\numberwithin{equation}{section}
\usepackage{amsthm}
\newtheorem{thm}{Theorem}[section]

\newtheorem{conj}[thm]{Conjecture}
\newtheorem{dfn}[thm]{Definition}

\usepackage[a4paper, top = 1.2in, bottom = 1.2in, left = 1.2in, right = 1.2in]{geometry}

\usepackage{enumitem}
\newlist{steps}{enumerate}{1}
\setlist[steps, 1]{label = Step \arabic*:}

\numberwithin{equation}{section}

\newcommand{\A}{\mathbb{A}}
\newcommand{\C}{\mathbb{C}}
\newcommand{\F}{\mathbb{F}}
\newcommand{\N}{\mathbb{N}}
\newcommand{\Q}{\mathbb{Q}}

\newcommand{\R}{\mathbb{R}}

\newcommand{\Z}{\mathbb{Z}}

\newcommand{\mbH}{\mathbb{H}}
\newcommand{\f}{\mathbf{f}}

\newcommand{\mcO}{\mathcal{O}}

\newcommand{\mfa}{\mathfrak{a}}
\newcommand{\mfb}{\mathfrak{b}}

\newcommand{\mfm}{\mathfrak{m}}
\newcommand{\mfn}{\mathfrak{n}}

\newcommand{\Dif}{\mathfrak{D}}
\newcommand{\mfp}{\mathfrak{p}}
\newcommand{\mfq}{\mathfrak{q}}
\newcommand{\mfP}{\mathfrak{P}}

\newcommand{\GL}{\mathrm{GL}}

\newcommand{\SO}{\mathrm{SO}}

\newcommand{\Gal}{\mathrm{Gal}}

\newcommand{\new}{\mathrm{new}}

\def\1{1\!\!1}

\newcommand{\bsmat}[4]{\bigl[ \begin{smallmatrix} #1 & #2 \\ #3 & #4 \end{smallmatrix} \bigr]}
\def\dis{\displaystyle}

\title[A survey on Diophantine equations over $K$]{A survey on the generalized Fermat equation of various signatures over totally real fields}

\author[S. Sahoo]{Satyabrat Sahoo}
\address[S. Sahoo]{Yau Mathematical Sciences Center, Tsinghua University, 
	Beijing 100084, China.}
\email{satyabrat.sahoo.94@gmail.com}

\date{\today}

\begin{document}
	\maketitle
	\begin{abstract}
		Following the famous proof of Fermat's Last Theorem by Andrew Wiles using the modularity of elliptic curves over $\mathbb{Q}$, significant developments have been made in the study of Diophantine equations using the modularity method. This article presents a survey of numerous results on the solutions of the generalized Fermat equation of signatures $(p,p,p)$, $(p,p,2)$, $(p,p,3)$, and $(r,r,p)$ over totally real number fields using the modularity method.
	\end{abstract}
	
	\maketitle
	
	\section{Introduction}
	The study of solutions to the Diophantine equations is one of the most attractive areas in number theory. The most prominent case is the Fermat equation $x^p+y^p=z^p$. In \cite{W95}, Wiles proved Fermat's Last Theorem and established that the Fermat equation $x^n+y^n=z^n$ has no non-trivial primitive integer solutions for integers $n \geq 3$.  
	Consider the generalized Fermat equation,
	\begin{equation}
		\label{generalized Fermat eqn}
		Ax^p+By^q+Cz^r=0, \text{ where } A,B,C, p,q,r \in \Z \setminus \{0\}
	\end{equation}
	with $\gcd(A,B,C)=1$ and primes $ p,q,r \geq 2$ with $\frac{1}{p} +\frac{1}{q}+ \frac{1}{r} <1$. We say $(p,q,r)$ as the \textbf{signature} of the generalized Fermat equation~\eqref{generalized Fermat eqn}.
	The following conjecture is known for the solutions of the equation~\eqref{generalized Fermat eqn} (see \cite{DG95}). 
	\begin{conj}
		\label{DG conj}
		Fix  $ A,B,C \in \Z $ with $\gcd(A,B,C)=1$. Then, over all choices of prime exponents $ p,q,r$ with $\frac{1}{p} +\frac{1}{q}+ \frac{1}{r} <1$, the the generalized Fermat equation~\eqref{generalized Fermat eqn} has only finitely many non-trivial primitive integer solutions.
	\end{conj}
	Here, a solution $(a,b,c) \in \Z^3$ to the equation~\eqref{generalized Fermat eqn} is called non-trivial if $abc \neq 0$, and primitive if $\gcd(a,b,c)=1$. In \cite{DG95}, Darmon and Granville proved Conjecture~\ref{DG conj} for fixed $ p,q,r$. More precisely:
	\begin{thm}  \rm(\cite[Theorem 2]{DG95})
		For fixed integers $ A,B,C \in \Z \setminus \{0\}$ and fixed primes $ p,q,r \geq 2$ with $\frac{1}{p} +\frac{1}{q}+ \frac{1}{r} <1$, the equation~\eqref{generalized Fermat eqn} has only finitely many non-trivial primitive integer solutions.	
	\end{thm}
	Throughout this article, we denote $K$ as a totally real number field and $\mcO_K$ as the ring of integers of $K$. 
	This article presents a summary of various results on the solutions of the generalized Fermat equation~\eqref{generalized Fermat eqn} with signatures $(p,p,p)$, $(p,p,2)$, $(p,p,3)$ and $(r,r,p)$ over $K$ using the modularity method. For an elliptic curve $E$ defined over $K$ and a rational prime $p$, let $\bar{\rho}_{E,p} : G_K:=\Gal(\bar{K}/K) \rightarrow \mathrm{Aut}(E[p]) \simeq \GL_2(\F_p)$ be  the mod-$p$ Galois representation of $G_K$ induced by the action of $G_K$ on the $p$-torsion points  $E[p]$ of $E$. For any Hilbert modular newform $f$ over $K$ of parallel weight $k$, level $\mfn$ with coefficient field $\Q_f$ and for any non-zero prime ideal $\omega$ of $\mcO_{\Q_f}$, let $\bar{\rho}_{f, \omega}: G_K \rightarrow \GL_2(\F_\omega)$ be the residual Galois representation attached to $f, \omega$. 
	\subsection{Structure of the article}
	In \S\ref{section for p,p,p}, \S\ref{section for p,p,2}, \S\ref{section for p,p,3} and \S\ref{section for r,r,p}, we give the literature for the solutions of the generalized Fermat equation~\eqref{generalized Fermat eqn} with signatures $(p,p,p)$, $(p,p,2)$, $(p,p,3)$ and $(r,r,p)$, respectively over $K$. In \S\ref{section strategy of proof}, we give the strategy of proof in the modular method. In \S\ref{section for Hilbert mod forms}, we define Hilbert modular forms. Finally
	in \S\ref{section for modularity}, we present a summary of numerous results regarding the modularity of elliptic curves over $K$.
	
	\section{Signature $(p,p,p)$ case}
	\label{section for p,p,p}
	In this section, we discuss various results on the solutions of the generalized Fermat equation~\eqref{generalized Fermat eqn} of signature $(p,p,p)$, specifically $Ax^p+By^p+Cz^p=0$ over totally real number fields $K$.
	\subsection{Over $\Q$}
	For $A=B=C=1$, Wiles first proved Fermat's Last Theorem. More precisely:
	\begin{thm} \rm{~\cite[Theorem 0.5]{W95}}
		For every positive integer $n \geq3$, the Fermat equation $x^n+y^n+z^n=0$ has no non-trivial primitive integer solution.
	\end{thm}
	Let $L\in \{3, 5, 7, 11, 13, 17, 19, 23, 29, 53, 59 \}$ and $r \geq 0$ be an integer. In~\cite{S87}, Serre proved that the equation $x^p+y^p+L^rz^p=0$ of exponent $p\geq 11$ with $p \neq L$ has no non-trivial primitive integer solutions. Furthermore, Serre proved that if $L \geq 3$ is a prime number which is neither a Fermat nor a Mersenne prime, then there exists a constant $C_L>0$ (depending on $L$) such that for all primes $p>C_L$, the equation $x^p+y^p+L^rz^p=0$ with integer $r >0$ has no non-trivial primitive integer solutions. In~\cite{DM97}, Darmon and Merel proved that the equation $x^n+y^n+2z^n=0$ with $n \geq 3$ has no non-trivial primitive integer solutions. In~\cite{R97}, Ribet demonstrated that the equation $x^p+y^p+2^rz^p=0$ has no non-trivial primitive integer solutions for $2\leq r <p$.
	Finally, in \cite{K97}, Kraus examined the integer solutions of the generalized Fermat equation $Ax^p+By^p+Cz^p=0$ for various choices of non-zero integers $A,B,C$.
	
	\subsection{Over $K$}
	Let $K$ be a totally real number field and $\mcO_K$ be the ring of integers of $K$. Consider the Diophantine equation
	\begin{equation}
		\label{Ax^p +By^p +Cz^p=0}
		Ax^p +By^p +Cz^p=0, 
	\end{equation}  
	of prime exponent $p \geq 5$, $\text{ where } A,B,C\in \mcO_K \setminus \{0\}$.
	\begin{dfn}[Trivial solution]
		We say a solution $(a, b, c)\in  \mcO_K^3$ to the equation~\eqref{Ax^p +By^p +Cz^p=0} is trivial, if $abc=0$, otherwise non-trivial. Further, we call it primitive if $a\mcO_K+b\mcO_K+c\mcO_K=\mcO_K$
	\end{dfn} 
	
	\begin{dfn}
		We say the Diophantine equation $Ax^p +By^p +Cz^p=0$ (resp. $Ax^p +By^p +Cz^r=0$, $r\in \{2,3\}$) has no asymptotic solution in a subset $S \subseteq \mcO_K^3$, if there exists a constant $V:=V_{K,A,B,C}$ (depending on $K,A,B,C$) such that for primes $p >V$, the equation $Ax^p +By^p +Cz^p=0$ (resp. $Ax^p +By^p +Cz^r=0$, $r\in \{2,3\}$) has no non-trivial primitive solutions in $S$.
	\end{dfn}
	
	Let $P$ be the set of all non-zero prime ideals of $\mcO_K$. For any subset $S \subseteq P$, let $\mcO_{S}:=\{\alpha \in K : v_\mfP(\alpha)\geq 0 \text{ for all } \mfP \in P \setminus S\}$ be the ring of $S$-integers in $K$ and $\mcO_{S}^*:= \{\alpha \in K : v_\mfP(\alpha)= 0 \text{ for all } \mfP \in P \setminus S\}$. We call the elements of $\mcO_{S}^*$ as the $S$-units. Let $S_K:=\{ \mfP \in P : \mfP|2 \}$, $T_K:=\{ \mfP \in S_K : \f(\mfP,2)=1 \}$ and $U_K:=\{ \mfP \in S_K: (3, v_\mfP(2))=1 \}$, where $\f(\mfP,2)$ denotes the residual degree of the prime ideal $\mfP$ lying above $2$. We now state the Eichler-Shimura conjecture. 
	\begin{conj}[Eichler-Shimura]
		\label{ES conj}
		Let $f$ be a Hilbert modular newform defined over $K$ of parallel weight $2$, level $\mfn$ and with coefficient field $\Q_f= \Q$. Then, there exists an elliptic curve $E_f /K$ of conductor $\mfn$ having the same $L$-function as $f$.
	\end{conj}
	In~\cite{D04}, Darmon proved that if either 
	$[K: \Q] $ is odd or there exists some prime ideal $\mfq \in P$ such that $v_\mfq(\mfn) = 1$, then Conjecture~\ref{ES conj} holds over $K$. We write $(ES)$ for ``either $[K: \Q]$ is odd or $T_K \neq \varphi$ or Conjecture \ref{ES conj} holds for $K$."  
	
	In \cite{FS15}, Freitas and Siksek first studied the asymptotic solution of the Fermat equation $x^p+y^p+z^p=0$ over $K$. More precisely:
	\begin{thm}\rm{\cite[Theorem 3]{FS15}}
		\label{main result for x^p+y^p=2^rz^p}
		Let $K$ be a totally real number field satisfying $(ES)$. Suppose, for every solution $(\lambda, \mu)$ to the $S_K$-unit equation
		\begin{equation}
			\label{S_K-unit solution}
			\lambda+\mu=1, \ \lambda, \mu \in \mcO_{S_K}^\ast,
		\end{equation}
		either 
		\begin{enumerate}
			\item there exists some $\mfP \in T_K$ that satisfies
			$\max \left\{|v_\mfP(\lambda)|,|v_\mfP(\mu)| \right\}\leq 4v_\mfP(2)$; or 
			\item there exists some $\mfP \in U_K$ that satisfies both $\max \left\{|v_\mfP(\lambda)|,|v_\mfP(\mu)| \right\}\leq 4v_\mfP(2)$ and $v_\mfP(\lambda \mu) \equiv v_\mfP(2) \pmod 3$.
		\end{enumerate}
		Then the equation $x^p+y^p+z^p=0$ has no asymptotic solution in $\mcO_K^3$.
	\end{thm}
	Moreover, in \cite{FS15}, the authors also proved that the equation $x^p+y^p+z^p=0$ has no asymptotic solution in $\mcO_K^3$ with $K=\Q(\sqrt{d})$ for a subset of $d \geq 2$ with density $\frac{5}{6}$ among the set of square-free integers $d \geq 2$.
	In \cite{D16}, Deconick extended the work of~\cite{FS15} to the equation $Ax^p+By^p+Cz^p=0$, where $A,B,C \in \mcO_K$ with $ABC$ is odd (in the sense that $\mfP \nmid ABC$, for all prime $\mfP \in S_K$). In \cite{KS24a}, Kumar and Sahoo studied the asymptotic solution of the equation $x^p+2^ry^p+z^p=0$ over $K$, for $r \in \N$. Finally, in \cite{S24}, the author extended the result of \cite{KS24a} to the generalized Fermat equation $Ax^p+By^p+Cz^p=0$, where $A,B,C \in \mcO_K \setminus \{0 \}$ with $ABC$ is even (in the sense that $\mfP | ABC$, for some prime $\mfP \in S_K$) and also calculated the density of all square free integers $d\geq 2$ such that the equation $Ax^p+By^p+Cz^p=0$ has no asymptotic solution over $K=\Q(\sqrt{d})$. Similar to Theorem~\ref{main result for x^p+y^p=2^rz^p}, all these results depend on some explicit bounds on the solution of the $S$-unit equation.
	
	\subsection{Construction of Frey elliptic curve}
	For any non-trivial solution $(a, b, c)\in \mcO_K^3$ to the equation $Ax^p+By^p+Cz^p=0$, the Frey elliptic curve is given by
	\begin{equation}
		\label{Frey curve for x^2=By^p+Cz^p of Type I}
		E/K : Y^2= X(X-Aa^p)(X+Bb^p),
	\end{equation}
	with $c_4=2^4(A^2a^{2p}-BCb^pc^p), \Delta_E=2^4A^2B^2C^2(abc)^{2p}$ and $ j_E=2^{8} \frac{(A^2a^{2p}-BCb^pc^p)^3}{A^2B^2C^2(abc)^{2p}}$, where $\Delta_E$ denotes as the discriminant of $E$ and $j_E$ denotes as the $j$-invariant of $E$.
	
	\section{Signature $(p,p,2)$ case}
	\label{section for p,p,2}
	In this section, we discuss various results on the solutions of the generalized Fermat equation~\eqref{generalized Fermat eqn} of signature $(p,p,2)$, i.e., $Ax^p+By^p+Cz^2=0$ over totally real number fields $K$.
	\subsection{Over $\Q$}
	In~\cite{DM97}, Darmon and Merel first proved that the equation $x^n+y^n=z^2$ with $n \geq 3$ has no non-trivial primitive integer solutions. In~\cite{I03}, Ivorra studied the integer solutions of the equation $x^2=y^p+2^rz^p$ and $2x^2=y^p+2^rz^p$ of exponent $p \geq 5$ with $0 \leq r < p$. In ~\cite{S03}, Siksek established that the only non-trivial primitive integer solutions to the equation $x^2=y^p+2^rz^p$  of exponent $p \geq 5$ with $r \geq 2$ are when $r=3,\ x=\pm3,\ y=z=1$. In \cite{BS04}, Bennet and Skinner studied the integer solutions of the equation $Ax^2=By^p+Cz^p$ for various choices of non-zero integers $A,B,C $. Recently in \cite{C24}, Cazorla Garc\`ia studied the integer solutions of the equation $Ax^2=q^ky^{2n}+z^n$ with $2|z$, for any fixed positive intgers $A,q,k$ with $A$ is square-free and $q\geq 3$ is a prime.
	
	\subsection{Over $K$}
	In~\cite{IKO20}, I\c{s}ik, Kara, and Ozman first studied the asymptotic solution of the equation $x^p+y^p=z^2$, and proved that  $x^p+y^p=z^2$ has no asymptotic solution $(a,b,c) \in \mcO_K^3$ with $2|b$ when the narrow class number $h_K^+=1$ and $T_K :=\{ \mfP \in S_K : \f(\mfP,2)=1 \}\neq \varphi$. In ~\cite{M22}, Mocanu generalized the result of \cite{IKO20} by replacing the assumption $h_K^+=1$ in ~\cite{IKO20} with $2 \nmid h_K^+$.
	In ~\cite{KS24a}, Kumar and Sahoo relaxed the assumptions in ~\cite{IKO20} and proved that $x^p+y^p=z^2$ has no asymptotic solution $(a,b,c) \in \mcO_K^3$ with $2|ab$ without assuming the conditions $h_K^+=1$ and $T_K \neq \varphi$. 
	Finally in \cite{KS24b}, Kumar and Sahoo extended the work of \cite{M22} and studied the asymptotic solutions of the equation $x^2=By^p+Cz^p$ over $K$, where $B$ is an odd integer and $C$ is either an odd integer or $2^r$ for some $r \in \N$. In \cite{KS24b}, the authors also studied the asymptotic solution of the equation $2x^2=By^p+2^rz^p$ over $K$, where $B$ is an odd integer and $r \in \N$. 
	
	\subsection{Construction of Frey elliptic curve}
	\begin{itemize}
		\item  For any non-trivial solution $(a, b, c)\in \mcO_K^3$ to the equation $x^2=By^p+Cz^p$, the Frey elliptic curve $E:=E_{a,b,c}$ is given by
		\begin{equation}
			\label{Frey curve for x^2=By^p+Cz^p of Type I}
			E/K : Y^2 = X(X^2+2aX+Bb^p),
		\end{equation}
		with $c_4=2^4(Bb^p+4Cc^p),\ \Delta_E=2^{6}(B^2C)(b^2c)^{p}$ and $ j_E=2^{6} \frac{(Bb^p+4Cc^p)^3}{B^2C(b^2c)^{p}}$.
		\item For any non-trivial solution $(a, b, c)\in \mcO_K^3$ to the equation $2x^2=By^p+Cz^p$, the Frey curve $E:=E_{a,b,c}$ is given by
		\begin{equation}
			\label{Frey curve for 2x^2=By^p+Cz^p of Type I}
			E/K: Y^2 = X(X^2-4aX+2Bb^p),
		\end{equation}
		with $c_4=2^5(Bb^p+2^{r+2}c^p),\ \Delta_E=2^{9+r}B^2(b^2c)^{p}$ and $ j_E=2^{6-r} \frac{(Bb^p+2^{r+2}c^p)^3}{B^2(b^2c)^{p}}$.
	\end{itemize}

	\section{Signature $(p,p,3)$ case}
	\label{section for p,p,3}
	In this section, we discuss various results on the solutions of the generalized Fermat equation~\eqref{generalized Fermat eqn} of signature $(p,p,3)$, i.e., $Ax^p+By^p+Cz^3=0$ over totally real number fields.
	\subsection{Over $\Q$}
	In~\cite{DM97}, Darmon and Merel proved that the equation $x^n+y^n=z^3$ with $n \geq 3$ has no non-trivial primitive integer solutions. In \cite{BVY04}, Bennett, Vatsal, and Yadzdani studied the integer solutions of the generalized Fermat equation of signature $(n,n,3)$, i.e., $Ax^n+By^n=Cz^3$, for various choices of non-zero integers $A,B,C $.
	
	\subsection{Over $K$}
	In \cite{M22}, Mocanu first studied the asymptotic solution $(a,b,c) \in \mcO_K^3$ to the equation $x^p+y^p=z^3$ with $3|b$ over $K$. 
	Recently in \cite{KS24c}, Kumar and Sahoo extended the work of \cite{M22} to the generalized Fermat equation of signature $(p,p,3)$, i.e., $Ax^p+By^p+Cz^3=0$ over $K$, where $A,B,C \in \mcO_K \setminus \{0\}$, and studied the asymptotic solution of the equation $Ax^p+By^p+Cz^3=0$ over $K$. 
	
	\subsection{Construction of Frey elliptic curve}
	For any non-trivial solution $(a, b, c)\in \mcO_K^3$ to the equation $Ax^p+By^p=Cz^3$ , the Frey elliptic curve is given by
	\begin{equation}
		\label{Frey curve for x^2=By^p+Cz^p of Type I}
		E/K: Y^2+3CcXY+C^2Bb^pY = X^3,
	\end{equation}
	where $c_4=3^2C^3c(9Aa^p+Bb^p),\ \Delta_E=3^3AB^3C^8(ab^3)^p$ and $ j_E=3^{3} \frac{Cc^3(9Aa^p+Bb^p)^3}{AB^3(ab^3)^p}$.
	
	\section{Signature $(r,r,p)$ case}
	\label{section for r,r,p}
	In this section, we discuss various results on the solutions of the Fermat equation of signature $(r,r,p)$, i.e.,
	\begin{equation}
		\label{r,r,p over Z}
		x^r+y^r=dz^p,
	\end{equation}
	over $K$, where $r,p\geq 5$ are rational primes and $d \in \mcO_K \setminus \{0\}$.
	
	\subsection{Over $\Q$}
	Throughout this subsection, we fix $d \in \Z \setminus \{0\}$.
	For $r=5$, the integer solutions of the equation~\eqref{r,r,p over Z} were initially studied by Billerey in~\cite{B07} and subsequently in~\cite{BD10}, \cite{DF14}, \cite{BCDF19}. In \cite{F15}, Freitas established a multi-Frey family of elliptic curves to study the integer solutions of the equation~\eqref{r,r,p over Z} and proved that the equation~\eqref{r,r,p over Z} (with $r=7$, $d=3$) has no non-trivial primitive integer solutions for primes $p> (1+3^{18})^2$. In \cite{BCDF23b}, Billerey et al. proved that the equation $x^7+y^7=3z^n$ (resp. $x^7+y^7=z^n$) with integers $n \geq 2$ has non-trivial primitive integer solution $(a,b,c)$ (resp. $(a,b,c)$ with $2 |a+b$ or $7 |a+b$).  
	In \cite{BCDF23a}, Billerey et al. examined the integer solutions of the equation $x^{11}+y^{11}=z^n$ with integers $n \geq 2$ using Frey abelian varieties. The integer solutions of the equation \eqref{r,r,p over Z} for $r=13$ were first studied in \cite{DF13}. In \cite{BCDF19} (resp. \cite{BCDDF23}), the authors proved the equation \eqref{r,r,p over Z} with $r=13$, $d=3$, and $p \neq 7$ (resp. for all $p$) has no non-trivial primitive integer solution.  
	
	In \cite{FN24}, Freitas and Najman proved that for a set of primes $p$ with positive density and for fixed integers $r, d$ with $r,p \nmid d$, the equation~\eqref{r,r,p over Z} has no non-trivial primitive integer solution $(a,b,c)$ of the form $2 |a+b$ or $r |a+b$. Recently in \cite{KMO24}, Kara, Mocanu and \"Ozman studied the asymptotic integer solutions of the equation \eqref{r,r,p over Z} by assuming Conjecture~\ref{ES conj} and the weak Frey-Mazur conjecture  (see \cite[Conjecture 3.11]{KMO24}).
	
	\subsection{Over K}
	For any prime $r \geq 5$, let $\zeta_r$ be a primitive $r$th root of unity in $\C$ and $K^+:= K(\zeta_r+ \zeta_r^{-1})$.
	\begin{dfn}[Trivial solution]
		We say a solution $(a, b, c)\in  \mcO_{K^+}^3$ to the equation~\eqref{r,r,p over Z} is trivial, if $abc=0$, otherwise non-trivial. Further, we call it primitive if $a\mcO_{K^+}+b\mcO_{K^+}+c\mcO_{K^+}=\mcO_{K^+}$
	\end{dfn} 
	
	\begin{dfn}
		We say the Diophantine equation $x^r+y^r= dz^p$ has no asymptotic solution in a subset $S \subseteq \mcO_{K^+}^3$ if there exists a constant $V:=V_{K,r,d}$ (depending on $K,r,d$) such that for primes $p >V$, the equation $x^r+y^r= dz^p$ has no non-trivial primitive solutions in $S$.
	\end{dfn}
	In \cite{M23}, Mocanu first studied the asymptotic solution $(a,b,c) \in \mcO_{K^+}^3$ of the equation $x^r+y^r=z^p$ with $2|c$. In \cite{JS25}, Jha and Sahoo extended the work of \cite{M23} to the equation~\eqref{r,r,p over Z}, and studied the asymptotic solution $(a,b,c) \in \mcO_{K^+}^3$ of the equation $x^r+y^r=dz^p$ with $2|c$ (resp. $2 |a+b$). In \cite{JS25}, the authors also studied the asymptotic solution $(a,b,c) \in \mcO_{K}^3$ of the equation $x^5+y^5=dz^p$ with $2 \nmid c$.
	
	\subsection{Construction of Frey elliptic curve for $x^r+y^r=dz^p$ over $K^+$ }
	In this subsection, we recall the Frey elliptic curve associated to any non-trivial primitive solution to the Diophantine equation $x^r+y^r=dz^p$ (see \cite{F15} for more details).
	Let $L:=K(\zeta_r)$.
	For any prime $r \geq 5$, choose $$\phi_r(x,y):= \frac{x^r+y^r}{x+y}= \sum_{i=0}^{r-1}(-1)^ix^{r-1-i}y^i.$$ 
	The factorization of the polynomial $\phi_r(x,y)$ over $L$ is given by
	\begin{small}
		\begin{equation*}
			\label{facored over L}
			\phi_r(x,y)= \prod_{i=1}^{r-1}(x+\zeta_r^i y).
		\end{equation*}
	\end{small}
	For any $k$ with $0 \leq k \leq \frac{r-1}{2}$, choose 
	$$f_k(x,y):= (x+\zeta_r^k y) (x+\zeta_r^{-k} y)= x^2+ (\zeta_r^{k}+\zeta_r^{-k})xy+ y^2.$$ Then $f_k(x,y) \in K^{+}[x,y]$ for all $k$.
	The factorization of the polynomial $\phi_r(x,y)$ over $K^{+}$ is given by
	\begin{small}
		\begin{equation*}
			\label{facored over K^+}
			\phi_r(x,y)= \prod_{i=1}^{\frac{r-1}{2}} f_k(x,y).
		\end{equation*}
	\end{small}
	Since $r\geq 5$, we get $\frac{r-1}{2} \geq 2$. Fix three integers $k_1, k_2, k_3$ such that $0 \leq k_1 <k_2 <k_3 \leq \frac{r-1}{2}$.
	Let 
	\begin{small}
		\begin{equation*}
				\label{eqn for alpha, beta, gamma}
				\alpha:= \zeta_r^{k_3}+\zeta_r^{-k_3}-\zeta_r^{k_2}-\zeta_r^{-k_2}, \nonumber\\
				\beta:= \zeta_r^{k_1}+\zeta_r^{-k_1}-\zeta_r^{k_3}-\zeta_r^{-k_3},\\
				\gamma:= \zeta_r^{k_2}+\zeta_r^{-k_2}-\zeta_r^{k_1}-\zeta_r^{-k_1} \nonumber.
		\end{equation*}
	\end{small}
	Then $\alpha, \beta, \gamma \in K^{+}$ and $\alpha f_{k_1}+ \beta f_{k_2}+\gamma f_{k_3}=0$. Let $A(x,y):=\alpha f_{k_1} (x,y)$, $B(x,y):=\beta f_{k_2} (x,y)$, and $C(x,y):=\gamma f_{k_3} (x,y)$. Hence $A(x,y), B(x,y),C(x,y) \in K^{+}[x,y]$.
	%
	Let $(a,b,c) \in \mcO_{K^+}^3$ be a non-trivial primitive solution of equation \eqref{r,r,p over Z}. Let $A:=A(a,b)$, $B:=B(a,b)$ and $C:=C(a,b)$. Then $A,B,C \in K^+$ and $A+B+C=0$. So the Frey elliptic curve is given by
	\begin{equation}
		\label{Frey curve result1}
		E/K^+: Y^2=X(X-A)(X+B),
	\end{equation}
	where $\Delta_E=2^4(ABC)^2$, $c_4=2^4(AB+BC+CA)$ and $j_E=2^8 \frac{(AB+BC+CA)^3}{(ABC)^2}$. The Frey curve in equation~\eqref{Frey curve result1} depends on $p$ since $A,B,C$ depend on $p$.	
	
	\subsection{Construction of Frey elliptic curve for $x^5+y^5=dz^p$ over $K$ }
	For any non-trivial primitive solution $(a, b, c)\in \mcO_K^3$ to the Diophantine equation $x^5+y^5=dz^p$, the Frey elliptic curve is given by
	\begin{equation}
		\label{Frey curve result2}
		E/K : y^2= x^3-5(a^2+b^2)x^2+5\phi_5(a,b)x,
	\end{equation}
	where $\phi_5(a,b):=\frac{a^5+b^5}{a+b}, \ c_4=2^4 5 (5(a^2+b^2)^2-3\phi_5(a,b)),\\ \Delta_E=2^4 5^3 (a+b)^2(a^5+b^5)^2$ and $ j_E=2^{8} \frac{(5(a^2+b^2)^2-3\phi_5(a,b))^3}{(a+b)^2(a^5+b^5)^2}$ (see \cite[\S2]{B07} for the Frey curve). 
	
	\section{Strategy of the proof}
	\label{section strategy of proof}
	In this section, we provide a brief strategy of the proof in the modular method for solving the generalized Fermat equation over totally real fields $K$. The following are the key steps in the modular method:
	\begin{steps}
		\item For any non-trivial solution $(a, b, c)\in \mcO_K^3$ to the equation \eqref{generalized Fermat eqn} of signature $(p,p,p)$ or $(p,p,2)$ or $(p,p,3)$ or $(r,r,p)$, associate a suitable \textbf{Frey elliptic curve $E/K$}.
		\item Then prove the \textbf{modularity} of $E/K$ for primes $p \gg 0$, $E/K$ has \textbf{semi-stable reduction} (i.e., either good or multiplicative reduction) at all primes $\mfq|p$, and the mod-$p$ Galois representation $\bar{\rho}_{E,p}$ is \textbf{irreducible} for $p \gg 0$.
		
		\item Using \textbf{level lowering results} of $\bar{\rho}_{E,p}$, we have $\bar{\rho}_{E,p} \sim \bar{\rho}_{f,p}$, for some Hilbert modular newform $f$ of parallel weight $2$ with rational eigenvalues of lower level $\mfn_p$ which strictly divides the conductor of $E$.
		
		\item Finally to get a \textbf{contradiction}, prove that the finitely many Hilbert modular newforms of level $\mfn_p$ that occur in Step $3$ do not correspond to $\bar{\rho}_{E,p}$ which contradicts $\bar{\rho}_{E,p} \sim \bar{\rho}_{f,p}$.
	\end{steps}
	
	We now explain the proof of Theorem~\ref{main result for x^p+y^p=2^rz^p} briefly in terms of the above steps of the modular method. 
	\begin{itemize}
		\item For any non-trivial solution $(a,b,c)\in \mcO_K^3$ to the equation $x^{p}+y^{p}+z^p=0$, the Frey elliptic curve $E/K$ is given by 
		$$E/K : Y^2= X(X-a^p)(X+b^p).$$
		
		\item By \cite[Corollary 2.1]{FS15}, the Frey curve $E/K$ is modular  for $p \gg 0$. By \cite[Lemma 3.3]{FS15}, the Frey curve $E/K$ is minimal and semistable at all prime ideals $\mfq|p$. By \cite[Theorem 6]{FS15}, we get $\bar{\rho}_{E,p}$ is irreducible for $p \gg 0$.
		
		\item  Then using the level-loweing result, i.e., \cite[Theorem 7]{FS15}, we get $\bar{\rho}_{E,p} \sim \bar{\rho}_{f,p}$, for some Hilbert modular newform $f$ of parallel weight $2$ with lower level $\mfn_p$ which strictly divides the conductor of $E$.
		
		\item We now explain the step $4$ briefly. 	
		First, by using \cite[Lemma 3.7]{FS15}, we conclude that the Frey curve $E/K$ has potential multiplicative reduction at all primes $\mfP \in T_K$, i.e., $v_\mfP(j_E)<0$ and $p | \# \bar{\rho}_{E,p}(I_\mfP)$  for $p > 4v_\mfP (2)$.
		Next, by using \cite[Theorem 7]{FS15}, there exists a constant $V_K>0$ (depending on $K$) and an elliptic curve $E'$ defined over $K$ such that 
		$E^\prime/K$ has good reduction away from $S_{K}$, $E^\prime/K$ has full $2$-torsion points, $\bar{\rho}_{E,p} \sim \bar{\rho}_{E^\prime,p}$ for primes $p \geq V_K$.
		
		Then, by using an image of inertia result, i.e., \cite[Lemma 3.4]{FS15}, we conclude that $E'$ has potential multiplicative reduction at all primes $\mfP \in T_{K}$, i.e., $v_\mfP(j_{E^\prime}) <0$ for all primes $\mfP \in T_{K}$.
		Next, we express the $j$-invariant $j_{E'}$ of $E'$ in terms of the solution $(\lambda, \mu)$ to the $S_{K}$-unit equation~\eqref{S_K-unit solution} by using the above properties of $E'$. Finally by using the explicit bound $\max \left\{|v_\mfP(\lambda)|,|v_\mfP(\mu)| \right\}\leq 4v_\mfP(2)$, we get $v_\mfP(j_{E^\prime}) \geq0$ for some $\mfP \in T_{K}$, which is a contradiction (see \cite[\S5]{FS15} for more details). A similar calculation also holds for primes $\mfP \in U_{K}$.
	\end{itemize}
	
	\section{Hilbert modular forms}
	\label{section for Hilbert mod forms}
	In this section, we define Hilbert modular forms.
	Let $F$ be a totally real field of degree $n$ over $\Q$ and let $\mcO_F$ be the ring of integers of $F$. Let $k=(k_1,\dots, k_n)\in \Z_{>0}^n$, and let $\mbH:=\{ z \in \C : \text{im(z)} >0\}$ be the upper half-plane. Let $\Dif$ be the absolute different of $\mcO_F$. For any non-zero prime ideal $\mfp$ of $\mcO_F$, 
	let $F_\mfp$ be the completion of $F$ at $\mfp$. For any ideals $\mfa$ and $\mfb$ of $\mcO_F$, we put
	\[ K_\mfp(\mfa, \mfb):=\left\{\bsmat{a}{b}{c}{d}\in \GL_2(F_\mfp)\, : \, \begin{matrix} a\in \mcO_\mfp, & b\in \mfa_\mfp^{-1}\Dif_\mfp^{-1}, & \\ c\in \mfb_\mfp\Dif_\mfp, & d\in\mcO_\mfp, & |ad-bc|_\mfp=1\end{matrix}\right\},\]
	where $\mfm_\mfp=\mfm \mcO_\mfp$ for any ideal $\mfm$ of $\mcO_F$ and $\Dif_\mfp=\Dif \mcO_\mfp$. Then $K_\mfp(\mfa, \mfb)$ is a subgroup of $\GL_2(F_\mfp)$. 
	Furthermore, we put 
	\[ K_0(\mfa, \mfb)=\SO(2)^n\cdot\prod_{\mfp<\infty}K_\mfp(\mfa, \mfb), \quad  \text{and} \quad W(\mfa, \mfb)=\GL_2^+(\R)^nK_0(\mfa, \mfb).\]
	In particular, if $\mfa=\mcO_F$, we simply write $K_\mfp(\mfb)=K_\mfp(\mcO_F, \mfb)$ and $W(\mfb)=W(\mcO_F, \mfb)$.
	Let $\A_F$ be the ad\`{e}le ring of $F$. Let $h$ be the narrow class number of $F$, and let $\{t_\nu\}_{\nu=1}^h$ taken to be a complete set of representatives of the narrow class group of $F$, where the archimedean part $t_{\nu, \infty}$ of $t_\nu$ is $1$ for all $\nu$.
	Then, we have the following disjoint decomposition of $\GL_2(\A_F)$
	\begin{equation}\label{eqn:decomp}
		\GL_2(\A_F)=\cup_{\nu=1}^h\GL_2(F)x_\nu^{-\iota} W(\mfb),
	\end{equation}
	where $\dis x_\nu^{-\iota} =\bsmat{t_\nu^{-1}}{}{}{1}$. 
	For each $\nu$, we define a congruence subgroup of $\GL_2(F)$ as follows:
	\begin{align*}
		\Gamma_\nu(\mfa, \mfb) &:= \GL_2(F)\cap x_\nu W(\mfa, \mfb)x_\nu^{-1} \\ 
		&= \left\{ \bsmat{a}{t_\nu^{-1}b}{t_\nu c}{d}\in\GL_2(F): \, \begin{matrix} a\in \mcO_F, & b\in\mfa^{-1}\Dif_F^{-1}, & \\ c\in \mfb\Dif_F, & d\in\mcO_F, & ad-bc\in \mcO_F^\times \end{matrix}\right\}. 
	\end{align*} 
	In particular, if $\mfa=\mcO_F$, we write $\Gamma_\nu(\mfb)=\Gamma_\nu(\mcO_F, \mfb)$ as before. 
	
	Let $z=(z_1, \dots, z_n) \in \mbH^n$. For a holomorphic function $f$ on $\mbH^n$ and an element $\gamma=(\gamma_1, \dots, \gamma_n) \in \GL_2(\R)^n$, define 
	$$  f ||_k \gamma(z):= \prod_{i=1}^{n} \det(\gamma_i)^{\frac{k_i}{2}} j(\gamma_i, z_i)^{-k_i} f(\gamma z),$$
	where $\gamma z=(\gamma_1z_1, \dots, \gamma_n z_n)$, $j(\gamma_i, z_i)=c_iz_i+d_i$ and $\gamma_iz_i= \frac{a_iz_i+b_i}{c_iz_i+d_i}$ for $\gamma_i=\left(\begin{matrix} a_i & b_i \\ c_i & d_i \end{matrix} \right)$.
	
	Let $M_k(\Gamma_\nu(\mfb))$ denote the space of all functions $f_\nu$ that are holomorphic on $\mbH^n$ as well as at cusps and satisfying
	\[ f_\nu ||_k \gamma=f_\nu \]
	for all $\gamma$ in $\Gamma_\nu(\mfb)$. We note that such a function $f_\nu$ has a Fourier expansion
	\begin{equation}
		\label{Fourier expansion of f_nu}
		f_\nu(z)=\sum_{\xi\in F}a_\nu(\xi) \exp(2\pi i \xi z)
	\end{equation}
	where $\xi$ runs over all the totally positive elements in $t_\nu^{-1}\mcO_F$ and $\xi=0$. We call $M_k(\Gamma_\nu(\mfb))$ as the space of all Hilbert modular forms of weight $k$ and level $\Gamma_\nu(\mfb)$.
	A Hilbert modular form $f_\nu \in M_k(\Gamma_\nu(\mfb))$ is a cusp form if the constant term of $f_\nu||_k\gamma$
	in its Fourier expansion is $0$ for all $\gamma \in \GL^+_2(F)$, and the space of such cusp forms $f_\nu$ with respect to $\Gamma_{\nu}(\mfb)$ is denoted by $S_k(\Gamma_\nu(\mfb))$.

	Let $M_k(\mfb)$ denote the space of all functions $f:=(f_1,\dots,f_h)$ where $f_\nu \in M_k(\Gamma_\nu(\mfb))$ for each $\nu\in \{1,\dots,h\}$ and satisfying
	\[ f(g)=f(\gamma x_\nu^{-\iota}w)=(f_\nu||_k w_\infty)(\mathbf{i})\]
	where $g=\gamma x_\nu^{-\iota}w\in\GL_2(F)x_\nu^{-\iota}W(\mfb)$ as in equation~\eqref{eqn:decomp} and $\mathbf{i}=(i,\dots, i)$. 
	We call $M_k(\mfb)$ the space of all Hilbert modular forms of weight $k$ and level $\mfb$. In particular, if $f_\nu \in S_k(\Gamma_\nu(\mfb))$ for all $\nu$, then the space of such $f$ is denoted by $S_k(\mfb)$ and called the space of all Hilbert modular cusp forms of weight $k$ and level $\mfb$.
	
	Let $\sigma_1, \dots \sigma_n$ be the embeddings of $F$ into $\R$. For any $m:= (m_1, \dots, m_n ) \in \Z^n$ and $\alpha \in F$, we write $\alpha^m=\prod_{i=1}^{n}(\sigma_i(\alpha))^{m_i}$.
	Let $\mfm$ be an ideal of $\mcO_F$ and write $\mfm=\xi t_\nu^{-1}\mcO_F$ with a totally positive element $\xi$ in $F$. We define the Fourier coefficients of $f$ as follows:
	\begin{equation}\label{eqn:coeff}
		C(\mfm,f) :=\begin{cases} N_{F/\Q}(\mfm)^{k_0/2}\xi^{-k/2}a_\nu(\xi) \quad & \text{if}\quad \mfm=\xi t_\nu^{-1}\mcO_F\subset \mcO_F, \\
			0 & \text{if} \quad \mfm \, \text{is not integral} \end{cases} 
	\end{equation}
	where $k_0=\max \{k_1,\dots,k_n\}$ and $a_\nu(\xi)$ as in equation~\eqref{Fourier expansion of f_nu}.
	
	We denote $S_k^\new(\mfn)$ for the space of all Hilbert modular newforms over $F$ of weight $k$ and level $\mfn$. We say $f \in S_k^\new(\mfn)$ is normalized if $C(\mcO_F,f)=1$.  A primitive Hilbert modular newform is a normalized Hecke eigenform in $S_k^\new(\mfn)$ (see \cite{S78} for more details).

	\section{Modularity of elliptic curves over $K$}    
	\label{section for modularity} 
	In this section, we recall various results on the modularity of elliptic curves over totally real number fields. We first recall the definition of the modularity of an elliptic curve. To do this, we first review the definitions of the L-function attached to an elliptic curve and a Hilbert modular form.
	\begin{dfn}[$L$-function attached to an elliptic curve]
		Let $K$ be a number field and $E/K$ be an elliptic curve. Then, the $L$-function attached to the elliptic curve $E$ is defined by the Euler product
		$$ L(E,s)= \prod_{\mfp} L_\mfp \left(	N_{K/ \Q}(\mfp)^{-s} \right)$$
		for $s \in \C$, where $\mfp$ runs over all prime ideals of $\mcO_K$,
		$$
		L_\mfp(T)^{-1} =  \begin{cases}
			1-a_{\mfp}T+ N_{K/ \Q}(\mfp)T^2 & \text{if $E$ has good reduction at $\mfp$,}\\
			1-a_{\mfp}T & \text{if $E$ has bad reduction at $\mfp$},
		\end{cases}$$
		and 
		$$
		a_{\mfp} =  \begin{cases}
			N_{K/ \Q}(\mfp)+1-\#E(\F_\mfp) & \text{if $E$ has good reduction at $\mfp$,}\\
			0 & \text{if $E$ has additive reduction at $\mfp$,}\\
			1 &  \text{if $E$ has split multiplicative reduction at $\mfp$,}\\
			-1 &  \text{if $E$ has non-split multiplicative reduction at $\mfp$.}
		\end{cases}$$
		(see \cite[Appendix C, \S16]{S86} for more details).
	\end{dfn}
	
	\begin{dfn}[$L$-function attached to a Hilbert modular form]
		Let $F$ be a totally real number field.
		Let $f \in S_k^\new(\mfn) $ be a primitive Hilbert modular form over $F$, of weight $k$ and level $\mfn$ with Fourier coefficients $C(\mfm,f)$. Then, the $L$-function attached to $f$ is defined by 
		$$ L(f,s)= \sum_{\mfm} C(\mfm, f)	N_{K/ \Q}(\mfm)^{-s}$$
		for $s \in \C$, where $\mfm$ runs over all ideals of $\mcO_F$.
	\end{dfn}
	Then $L(f,s)$ has an Euler product given by
	$$ L(f,s)= \prod_{\mfp} \left(1-C(\mfp, f) N_{K/ \Q}(\mfp)^{-s}+N_{K/ \Q}(\mfp)^{k_0-1-2s} \right)^{-1},$$
	where $\mfp$ runs over all prime ideals of $\mcO_F$ (see \cite[\S5.3]{T12} for more details).
	
	\begin{dfn}[Modularity of elliptic curve]
		Let $K$ be a totally real number field. We say an elliptic curve $E/K$ is modular over $K$ if there exists a primitive Hilbert modular newform $f$ over $K$ of parallel weight $2$ and level $\mfn$ with rational eigenvalues such that both $E$ and $f$ have same $L$- function, i.e., $  L(E,s)=L(f,s)$.
	\end{dfn}
	
	The Shimura-Taniyama conjecture states that every elliptic curve over $\Q$ is modular. This conjecture was first proved by Wiles \cite{W95}, and Taylor and Wiles \cite{TW95} for semi-stable elliptic curves over $\Q$. Later, in \cite{BCDT01}, Breuil, Conrad, Diamond, and Taylor proved this conjecture for all elliptic curves over $\Q$. 
	\begin{thm}
		Every elliptic curve $E$ over $\Q$ is modular.
	\end{thm}
	The following results are known for the modularity of elliptic curves over $K$ when $[K: \Q] \in \{ 2,3,4,5 \}$. More precisely:
	\begin{thm}\rm(\cite[Theorem 1]{FLHS15})
		Every elliptic curve $E$ over a real quadratic field is modular.
	\end{thm}
	\begin{thm}\rm(\cite[Theorem 2]{DNS20})
		Every elliptic curve $E$ over a real cubic field is modular.
	\end{thm}
	\begin{thm}\rm{\cite[Theorem 1.1]{B22}}
		Let $K$ be a totally real quartic field not containing $\sqrt{5}$. Then, every elliptic curve $E$ over $K$ is modular.
	\end{thm}
	\begin{thm}\rm{\cite[Theorem 1.2]{IIY22}}
		For all but finitely many totally real fields $K$ of degree $5$ over $\Q$, every elliptic curve $E$ over $K$ is modular.
	\end{thm}
	In general, for any totally real number field $K$, \cite{FLHS15} proved that all but finitely many elliptic curves over $K$ are modular. More precisely:
	\begin{thm}\rm{\cite[Theorem 5]{FLHS15}}
		Let $K$ be a totally real number field. Then, all but finitely many $\bar{K}$-isomorphism classes of elliptic curves over $K$  are modular. 
	\end{thm}
	\section*{Acknowledgments} 
	The author expresses his sincere gratitude to Prof. Narasimha Kumar and Prof. Nuno Freitas for the invaluable assistance in understanding the article~\cite{FS15}. The author sincerely thanks the reviewer for the insightful comments and valuable suggestions towards the improvement of this article.

\end{document}